\documentclass[12pt]{article}
\usepackage{amsmath}
\usepackage{amsfonts}
\usepackage{amssymb}
\usepackage[latin1]{inputenc}
\usepackage{latexsym}
\usepackage{amscd}
\usepackage{graphicx}

\input amssym.def
\newcommand {\C} {\mathbb{C}}
\newcommand {\N} {\mathbb{N}}

\newcommand {\R} {\mathbb{R}}

\newcommand {\sobre} {\longrightarrow}
\newcommand {\tende} {\rightarrow}

\newtheorem{df}{Definition}[section]

\newtheorem{prop}[df]{Proposition}
\newtheorem{teo}[df]{Theorem}
\newtheorem{corol}[df]{Corollary}
\newtheorem{lemma}[df]{Lemma}

\newtheorem{rema}[df]{Remark}

\begin{document}

\title{\textsc{New classes of spaceable sets of analytic functions on the open unit disk}}
\author{M. L. Louren\c co and D. M. Vieira\\USP, S\~ao Paulo, Brazil}
\date{}
\maketitle

\abstract{ In this paper we study an algebraic and topological structure  inside the following sets of special functions:  Bloch functions defined on the open unit disk that are unbounded and  analytic functions of bounded type  defined  a Banach algebra $E$ into  $E$, which are not Lorch-analytic.  }

\section{Introduction}

The concept of linearity was probably first idealized by Gurariy in \cite{gur1} (see also \cite{gur2, arogursep}), where he shows that the set of continuous functions on the closed unit interval [0,1] that are not differentiable at any point (\emph{continuous everywhere, differentiable nowhere functions}) contains, except for the null function, an infinite dimensional vector space. The first example of such a function was given by K. Weierstrass in 1872, and the function became known as \emph{Weierstrass' Monster}. Gurariy's result is important because it shows that there are lots of such "monsters", which could be thought as the opposite, that is, that these pathological functions could be rare. And after that, several results of this type were obtained, and this subject is still a field of big interest. Actually, in the last two decades we have seen a crescent interest in the search of linear structures within sets (mainly sets of functions or sequences) that do not enjoy themselves such structures. Other similar concepts naturally appeared, such as spaceability, maximal lineability, algebrability and others. As until 2016 there was a large number of paper published on the subject in the literature, R. M. Aron, L. Bernal-Gonz\' alez, D. M. Pellegrino and J. B. Seoane-Sep\' ulveda in  \cite{aroberpelsep} published a book to report the progress in the theme so far. After that, the subject continues in studies by many authors and in different directions.

In this paper we wish to enlarge the list of examples in some functions spaces. Indeed, we will show results on spaceability and residuality on some sets of analytic functions related to Bloch functions, Bergman functions and Lorch analytic mappings. To be more precise, let us fix the notation. Let $D$ be the open unit disk on the complex plane. The Banach space of all \emph{Bloch functions} will be denoted by $\mathcal{B}$, and the Banach space of \emph{little Bloch functions} will be denoted by $\mathcal{B}_0$, both endowed with the Bloch norm. The Banach space of all bounded analytic functions $f:D\sobre\C$ will be denoted by $\mathcal{H}^{\infty}(D)$, endowed with the \emph{sup} norm. It is known that $\mathcal{H}^{\infty}(D)\subset \mathcal{B}$, but the inclusion is strict, and the set $(\mathcal{B}\setminus \mathcal{H}^{\infty}(D)) \cup \{0\}$ is not a vector space  in $\mathcal{B}.$ The same is true about the inclusion $\mathcal{B}_0\subset\mathcal{B}$. We investigate the sets $\mathcal{B}\setminus \mathcal{H}^{\infty}(D)$ and $\mathcal{B}\setminus \mathcal{B}_0$, and we show that they are spaceable and residual in $\mathcal{B}$. We refer to \cite{andclupom, cima, tim, zhu1, zhu2} for background on Bloch and little Bloch functions.

In the last section we study analytic functions defined on Banach algebras. For a complex commutative Banach algebra $E$, the Fr\'echet space of all holomorphic functions of bounded type from $E$ into $E$ will be denoted by $\mathcal{H}_b(E,E)$, endowed with the topology of the uniform convergence on the bounded subsets of $E$. We denote the set of all Lorch-analytic mappings from $E$ into $E$ by $\mathcal{H}_L(E,E)$. The class of Lorch-analytic mappings (cf. Definition 2.1) was introduced by E. R. Lorch in \cite{lor}. It was shown in \cite{luizalex} that $\mathcal{H}_L(E,E)$ is a closed subalgebra of $\mathcal{H}_b(E,E)$. However, this inclusion is strict, even in the two-dimensional case. The authors in \cite{lourvi2}  proved that for $E=\C^2$, the set $\mathcal{H}(\C^2,\C^2)\setminus\mathcal{H}_L(\C^2,\C^2) $ is not a vector space, yet is spaceable but not residual in $\mathcal{H}(\C^2,\C^2)$. In this note, Section 3, we investigate this set in the case where E is  a general Banach algebra $E$, and we show that it is spaceable, but not residual in $\mathcal{H}_b(E,E)$. For background on Lorch-analytic mappings, we  refer to \cite{lor, luizalex}.

Finally, let us recall some  necessary definitions, which nowadays have became usual terminology.
Assume that $Y$ is a complex  vector space. Then   a subset $A$ of $Y$ is called: \textbf{lineable} if $A\cup\{0\}$ contains an infinite dimensional vector space;  \textbf{maximal lineable} if $A\cup\{0\}$ contains a vector subspace $S$ of $Y$ with $dim(S) = dim(Y).$

If, in addition, $Y$ is a topological vector space, then a subset $A$ of $Y$  is called:  \textbf{spaceable } if $A\cup\{0\}$ contains  some infinite dimensional closed vector space. If $Y$ is a Fr\'echet space, then a subset $A$ of $Y$  is called \textbf{residual in} $Y$ if $Y\setminus A=\cup_{n=1}^{\infty}F_n$, with $\stackrel{\circ}{\overline{F_n}}=\emptyset$. So, by Baire's theorem, residual sets are topologically large. Other interesting property  is  algebrability. If $Y$ is a function algebra, $A\subset Y$ is said to be: \textbf{algebrable} if there is an algebra $\mathcal{B}\subset A\cup\{0\}$, such that $\mathcal{B}$ has an infinite minimal system of generators. For background on above concepts we refer to \cite{aroberpelsep, arosep, berboni, bercab}.

\section{Bloch spaces}

The space of Bloch functions has been studied by many authors, because of its intrinsic interest and because it is the meeting place of
several areas of analysis.  An important property of the Bloch space is that it is invariant
under Mobius transforms. Another interesting aspect of Bloch functions on the unit disk is that they form a Banach space which has several noteworthy properties. For example, it is isomorphic to the classical Banach space  $\ell_{\infty}.$

In this section we will also study the Bergman spaces, both for being related to the Bloch space and for its intrinsic value. The Bloch space can be naturally identified with the dual of some Bergman space. The Bergman spaces brings together complex function theory, functional analysis and operator theory. It comes in contact with harmonic analysis, approximate theory, hyperbolic geometry, potential theory and partial differential equation. Our interest is to analyze how big is the topological and algebraic structure between the spaces of Bloch and Bergman functions.

We are interested in the classical Bloch space. In classical function theory on the complex  open unit disk, the Bloch space is a central object of study. The examples of Bloch functions are the set of polynomials and also the bounded analytic functions. We suggest \cite{andclupom, cima, tim, zhu1, zhu2} to see the basic concepts and properties of Bloch and little Bloch functions.

Let $D$ denote the open unit disk in the complex plane, let $\mathcal{H}(D)$ be the space of all analytic functions on $ D$  and $\mathcal{H}^{\infty}(D)$ be the subspace of $\mathcal{H}(D)$ of all bounded functions.  The space of \emph{Bloch functions} is defined in the following way
$$\mathcal{B}=  \Big\{ f \in \mathcal{H}(D) : \sup_{ \vert z\vert < 1}(1- \vert z\vert^2) \vert f^{'}(z)\vert < \infty\Big\},$$
and the Bloch norm is $||f||_{\mathcal{B}} = \vert f(0)\vert +   \sup_{ \vert z \vert < 1}(1- \vert z\vert^2)\vert f'(z)\vert$. 

We call $\mathcal{B}$, endowed with the norm $\| \cdot  \|_{\mathcal{B}}$, the Bloch space. It is well-known that $\mathcal{B}$ is a (nonseparable) Banach space. A Bloch function $f$ is an analytic function on $D$ whose derivative grows so faster than a constant times the reciprocal of the distance from $z$ to the boundary of $D$. 

The \emph{little Bloch space}, $\mathcal{B}_0$, is the closed separable subspace of $\mathcal{B}$  given by
 $$\mathcal{B}_0 = \{ f\in\mathcal{B}\,\,\,:\, \,\lim_{|z|\tende 1-}(1-|z|^2)|f'(z)|=0\}.$$ 

The space $\mathcal{B}_0$  has the same relation to $\mathcal{B},$ roughly speaking, as the disk algebra, of all analytic functions in $D$ and continuous in $\overline{D}$, has to  $\mathcal{H}^{\infty}(D).$  However $\mathcal{B}_0$  is not a Banach algebra under pointwise multiplication. 

There is a way to display several functions in the spaces $\mathcal{B}_0$ and $\mathcal{B}$, this can be done through lacunary series. A sequence $(\beta_n)_n$ of positive integers is called  a {\it gap sequence}  if there exists a constant $\beta > 1$ such that $\frac{\beta_{n+1}}{\beta_n} \geq \beta$ for all $n \in \N$.  In this case, we call a power series of the form  $\sum_{n= 0}^{+\infty} a_n z^{\beta_n}$ a {\it lacunary series,} where $(a_n)_n$ is a sequence of complex numbers. It is known that a lacunary series defines a function in $\mathcal{B}$ if and only if the coefficients $(a_n)$ are bounded. Similarly, a lacunary series defines a function in $ \mathcal{B}_0$ if and only if the coefficients $(a_n)$ tend to $0$ \cite[Theorem 3.15]{zhu1}. Indeed, the first known, and perhaps most striking, examples of functions in 
$\mathcal{B}$, which are not in any of the classical spaces are provided by gap power series. 

It is well-known the stric inclusion of $\mathcal{H}^{\infty}(D)$ in $\mathcal{B}$. The function $g:D \sobre\C $ defined by $g(z)=\log(1-z)$, for all $z\in D$ is an unbounded Bloch function. It is known that the functions $f_t(z)=\dfrac{e^{-it}}{2}$log$\left(\dfrac{1+e^{-itz}}{1-e^{-itz}}\right)$, $z\in D$, $0\leq t \le 2\pi$, are also unbounded Bloch functions. On the other hand, there is no inclusion relation between $\mathcal{H}^{\infty}(D)$ and $\mathcal{B}_0$. An example of a bounded analytic function on $D$ that does not belong to $\mathcal{B}_0$ can be found in \cite[Section 3]{atele}. The example of an unbounded function in $\mathcal{B}_0$ will be discussed after Theorem \ref{properties}. The inclusion $\mathcal{B}_0\subset\mathcal{B}$ is also strict. The function $f(z)=\log(1-z)$ is an element of the set $\mathcal{B}\setminus\mathcal{B}_0$.

Our aim is to study the algebraic and topological structures of the following two sets: $\mathcal{F}=\mathcal{B}\setminus\mathcal{H}^{\infty}(D)$ and $\mathcal{T}=\mathcal{B}\setminus\mathcal{B}_0$. We will show that $\mathcal{F}$ and $\mathcal{T}$ are spaceable and residual in $\mathcal{B}$. To obtain that $\mathcal{F}$ is spaceable, practically contradicts  the  following assertion  in \cite{andclupom} which says that "$\mathcal{B}$ is a little bigger than the space $\mathcal{H}^{\infty}(D)$". In the linear and topological sense, it could be said that $\mathcal{B}$ is \emph{quite} bigger than the space $\mathcal{H}^{\infty}(D).$

For the convenience of the reader, we include the next theorem and the following definition: let $f:D\sobre\C$ be an analytic function. For each $\alpha\in\C$, $|\alpha|\leq 1$, we define the  function $f_{\alpha}:D\sobre \C$ by $f_{\alpha}(z)=f(\alpha z)$, for all $z\in D$. When $r\in\R$, $0<  r<1$, $f_r$ is called the \emph{dilate function of} $f$. 

\begin{teo}\label{properties} 
\begin{enumerate}
\item $\mathcal{H}^{\infty}(D)\subset\mathcal{B}$ and $\|f\|_{\mathcal{B}}\leq\|f\|_{\infty}$;
\item  $\max_{|z|\leq r}|f(z)|\leq M_r \|f\|_{\mathcal{B}}$, for all $f\in\mathcal{B}$, $0<r<1$, where $M_r=1+\dfrac{1}{2}$log$\left(\dfrac{1+r}{1-r}\right)$.
\item  $f\in\mathcal{B}_0$ if, and only if, $\|f_r-f\|_{\mathcal{B}}\tende 0$, $r \tende 1^-$, where $f_r$ is the dilate function of $f$.

\end{enumerate}
\end{teo}

We observe that the proof of the  statements in Theorem \ref{properties} can be found in \cite{andclupom, zhu2}.

In order to discuss spaceability of $\mathcal{F}$, we observe that the space $\mathcal{H^{\infty}}(D)$, when endowed with the Bloch norm is not a Banach space. This happens because the classical Dirichlet space $\mathcal{D}\subset\mathcal{B}_0,$ so if $\varphi$ is the Riemannian function from $D$ onto a simply connected unbounded region in $\C$ with finite area, then $\varphi$ is a unbounded function in $\mathcal{B}_0$. Let $0<\alpha_n<1$ be such that $\alpha_n\tende 1$. Then $(\varphi_{\alpha_n})_{n\in\N}\subset\mathcal{H}^{\infty}(D)$ and by Theorem \ref{properties}, we have that $\varphi_{\alpha_n}\tende\varphi$ in the Bloch norm, but $\varphi \notin \mathcal{H}^{\infty}(D)$.

\begin{prop}\label{spaceable} $\mathcal{F}$ is spaceable.\end{prop}

\textbf{Proof:} Since $(\mathcal{B}, \Vert  \cdot \Vert_{\mathcal{B}})$ and $(\mathcal{H^{\infty}}(D), \Vert \cdot \Vert _{\infty})$ are Banach spaces  and $(\mathcal{H^{\infty}}(D), \Vert \cdot \Vert _{\mathcal{B}})$ is not closed in $\mathcal{B},$ the inclusion $(\mathcal{H^{\infty}}(D), \Vert \cdot \Vert_{\infty}) \hookrightarrow (\mathcal{B}, \Vert  \cdot \Vert_{\mathcal{B}})$ is continuous  by Theorem \ref{properties}. Since its image  is not closed, by applying \cite[Proposition 2.4]{kittim} we have that $\mathcal{F}$ is spaceable. $\Box$

 Naturally, if $\mathcal{F}$ is spaceable then it implies $\mathcal{F}$ is lineable, but here we decide to include the next result because it gives a construction of  linearly independent set in $\mathcal{F}$, from the  function $g(z)=\log(1-z)$, $z\in D$.

\begin{prop}\label{proprinc} There is an infinite linearly independent set  $S$  in $\mathcal{B}$ with $[S]\subset\mathcal{F}\cup\{0\}$.\end{prop}

\textbf{Proof:} Let $g:D \sobre\C $ be given by $g(z)=\log(1-z)$, for all $z\in D$, and, for each $\alpha \in \C$ with $ \vert \alpha \vert = 1$, we define the function by $g_{\alpha}(z) = g(\alpha z)$, for all $z\in D$. Consider $S= \{g_{\alpha}, \vert \alpha \vert = 1 \}.$ Consider now the vector space generated by $S$, denoted by $[S]$. We will prove that $S$ is a  linearly independent set and $[S]\subset\mathcal{F}\cup\{0\}$.  

Suppose that $f(z) = \sum_{k=1}^{N} c_k g_{\alpha_k} = 0$, $c_1,\dots,c_N\in\C$ and $\alpha_1,\dots,\alpha_N\in\C$ pairwise distinct with modulus $1$. In particular $f^{(k)}(0)=0$ (the $k$-th derivative of $f$ at $0$), for $k=1,\cdots,N$. It follows that $\sum_{k=1}^N c_k \alpha_k^j=0$, for $j=1,\dots N$, and by the Vandermonde determinant, we get that $c_1=c_2=\cdots=c_N$, and then $S$ is linearly independent.

If $f\in [S]$, then $f = \sum_{k=1}^N  c_kg_{\alpha_k},$ where $c_k \in \C, k = 1, \dots, N$ and $\alpha_1,\dots,\alpha_N\in\C$ pairwise distinct with modulus $1$. So  
$$\vert f'(z) \vert ( 1- \vert z \vert ^{2}) \leq \sum_{k=1}^N \frac{\vert c_k  \vert}{\vert 1 -\alpha_k z\vert}(1- \vert z \vert ^2)  \leq  \sum_{k=1}^N  \vert c_k \vert (1+\vert z \vert) \leq 2  \sum_{k=1}^N \vert c_k \vert.$$

This implies $f\in  \mathcal{B}.$ On the other hand, $ f \notin \mathcal{H^{\infty}}(D),$ because, if we take  $w_n  = \frac{1}{\alpha_1} ( 1 - \frac{1}{n})$  we have, for $k\neq 1$, that $1-\alpha_k w_n$ goes to $1- \frac{\alpha_k}{\alpha_1} \neq 0 $  and $|f(w_n)|$ diverges when  $n\tende\infty.$ Consequently $[S] \subset \mathcal{F}.$ $\Box$

As a consequence of Propostion \ref{proprinc} and the fact that dim $\mathcal{B} = \mathfrak{c}$, we have that $\mathcal{F}$ is maximal lineable. Next we show that $\mathcal{F}$ is topologically big.

\begin{prop} $\mathcal{F}$ is residual in $\mathcal{B}$.\end{prop}

\textbf{Proof:} Let $S_n = \{ f \in \mathcal{B}: \exists z \in D \textrm{ with } \vert f(z)\vert > n \}.$  Clearly, $\mathcal{F} = \cap_{n=1}^{\infty} S_n$ and it is sufficient to show that $S_n$ is an open set which is dense in $\mathcal{B}$.  Let $f \in S_n.$ Then there are $z \in D$, $\delta > 0 $ such that $\vert f(z)\vert = n +\delta$ and $0<r<1$ such that $|z|<r$. Consider  $0< \epsilon < \frac{\delta}{M_r}$ and $ g \in \mathcal{F}$ with $\Vert g- f \Vert_{\mathcal{B}} < \epsilon.$ Thus $|g(z)|>n$,  $g\in S_n,$ and as a consequence $S_n$ is open.

Let $h \in \mathcal{B}$ and $\varepsilon>0$ be given. If $h \notin \mathcal{H^{\infty}}(D),$ so $h \in S_n$, for all $n\in\N$. If $h \in \mathcal{H^{\infty}}(D),$ consider $\tilde{f} = \frac{\varepsilon}{3} \log (1-z)$. Then $ \tilde{f} \in \mathcal{F}$ and $\Vert \tilde{f} \Vert_{\mathcal{B}} = \frac{ 2\varepsilon}{3}.$ Let $\tilde{g} = \tilde{f}+h$ so $\tilde{g}\in S_n$ and $\Vert h-\tilde{g}\Vert _{\mathcal{B}} < \varepsilon$. $\Box$

\medskip

The space $\mathcal{B}$ is not an algebra under the usual multiplication, because   for example, the square of the Bloch function $\log(1+z)$ does not belong to $\mathcal{B}.$ A function $f$ is called a \emph{pointwise multiplier of} $\mathcal{B}$ if  for every $g \in \mathcal{B}$  the pointwise product $fg$ also belongs to  $\mathcal{B}.$ Thus  we denote a pointwise multiplier $f$ of the space $ \mathcal{B}$ by $f\mathcal{B} \subset \mathcal{B}.$ In \cite{zhu1} there is  a characterization of the pointwise multipliers of the Bloch space. We include it for the reader's convenience.

\begin{teo}\label{pointwise}\cite[Theorem 3.21]{zhu1} For an analytic function $f\in\mathcal{B}$, we have that $f\mathcal{B} \subset \mathcal{B}$ if, and only if, $f \in \mathcal{H}^{\infty}(D)$ and the function $(1- \vert z \vert^2) \vert f'(z)\vert \log \dfrac{1}{1- \vert z \vert^2}$ is bounded in $D.$\end{teo}

In particular, if $f$ is a pointwise multipier of $\mathcal{B}$, then $f$ must be bounded in $D$. In the following result, we present a consequence of Theorem \ref{pointwise}.

\begin{prop}\label{notalg} If $\mathcal{A}\subset\mathcal{B}$ is an algebra, then $\mathcal{A}\subset\mathcal{H}^{\infty}(D)$. \end{prop}

\textbf{Proof:} Since $\overline{\mathcal{A}}^{\|\cdot\|_{\mathcal{B}}}\subset\mathcal{B}$, we can assume that $(\mathcal{A},\|\cdot\|_{\mathcal{B}})$ is a closed algebra. Let $f\in \mathcal{A}$ and consider the linear operator $T:\mathcal{A}\sobre \mathcal{A}$ given by $T(g)=fg$, for all $g\in \mathcal{A}$. To show that $T$ is continuous, consider a sequence $(g_n)\subset\mathcal{A}$ such that $g_n\tende g\in \mathcal{A}$ and  suppose that $T(g_n)\tende h$. Using the  Theorem 2.1 (2), we  have that, the Bloch convergence implies pointwise convergence, that means, if we fix $z\in D$, and let $0<r<1$ be such that $|z|<r$, we have that $|g_n(z)-g(z)|\leq M_r\|g_n-g\|_{\mathcal{B}}$. Then $f(z)g_n(z)\tende f(z)g(z)$, and as $T(g_n)\tende h$, we have that  $f(z)g_n(z)\tende h(z)$. Consequently $h(z)=f(z)g(z)$ and by the Closed Graph Theorem, we have that $T$ is continuous.

Fix again $z\in D$, $0< r<1$ such that $|z|<r$ and consider the evaluation homomorphism $\delta_z:\mathcal{A}\sobre\C$ . Then $\delta_z\in (\mathcal{A},\|\cdot\|_{\mathcal{B}})'$. Consider now the adjoint $T^{\ast}:\mathcal{A}'\sobre \mathcal{A}'$. We have that $$T^{\ast}(\delta_z)(g)=(\delta_z\circ T)(g)=\delta_z(T(g))=\delta_z(fg)=f(z)g(z)=f(z)\delta_z(g),\,\forall g\in \mathcal{B}.$$ It follows $\|f(z)\delta_z\|=\|T^{\ast}(\delta_z)\|\leq\|T^{\ast}\|\|\delta_z\|$, which implies that $|f(z)|\leq\|T^{\ast}\|$, that is, $f\in \mathcal{H}^{\infty}(D)$.$\Box$

As we can see by the above proposition, the set $\mathcal{F}$ cannot be algebrable. We finish this section investigating the set $\mathcal{T}=\mathcal{B}\setminus\mathcal{B}_0$.

\begin{prop}\label{bzero} The set $\mathcal{T}=\mathcal{B}\setminus\mathcal{B}_0$ is spaceable.\end{prop}

\textbf{Proof:} For each $\alpha \in \C$ with $\vert \alpha \vert = 1,$ consider  $g_{\alpha}(z)=\log(1-\alpha z)$, for all $z\in D$. Let $S=\{g_{\alpha}\,:\,|\alpha|=1\}$.  By Propositon \ref{proprinc}, we have that $S$ is linearly independent and that $[S]\subset\mathcal{B}$. Here we  have to show that $[S]\subset\mathcal{T}$. If $f\in [S]$, then there exists $c_1,\dots,c_N\in\C$ and $\alpha_1,\dots,\alpha_N\in\C$ with modulus $1$ such that $f(z)= \sum_{k= 1}^{N}c_k\log(1-\alpha_k z).$ We can assume that $c_1\neq 0$ and that $\alpha_k$'s are distinct. Let $w_n=\frac{1}{\alpha_1}\left(1-\frac{1}{n}\right)$. Then $\lim_{n\tende\infty}(1-|w_n|^2)\left|\dfrac{c_1\alpha_1}{1-\alpha_1 w_n}\right|=2c_1$ and, for $k=2,\dots,N$, $\lim_{n\tende\infty}(1-|w_n|^2)\left|\dfrac{c_k\alpha_k}{1-\alpha_k w_n}\right|=0$. Then we have that $\lim_{n\tende\infty}(1-|z|^2)|f'(z)|\neq 0$, and hence $f\notin\mathcal{B}_0$. 

Now, we  consider the set of classes $\mathcal{C}=\{ \widehat{g_{\alpha}} : |\alpha|=1\}$ contained in the quotient space $\mathcal{B}/\mathcal{B}_0$. Suppose that $\sum_{k=1}^{N}\beta_k\widehat{g_{\alpha_k}}=\widehat{0}$, where $\alpha_k, \beta_k \in \C,$ $|\alpha_k|=1$, for $k= 1, \dots, N$. This implies that $f=\sum_{k=1}^{N}\beta_k g_{\alpha_k}\in\mathcal{B}_0$ and $f\equiv 0$, because if $f\neq 0$, then we would have that $f\in\mathcal{T}$, which is a contradiction. Since $S$ is linearly independent, we have that $\beta_1=\cdots=\beta_N=0$, and then the family $\mathcal{C}$ is linearly independent, so $\mathcal{B}/\mathcal{B}_0$ is infinite dimensional. Since $\mathcal{B}_0$ is closed in $\mathcal{B}$, it follows by \cite[Theorem 2.2]{kittim} that $\mathcal{T}$ is spaceable. $\Box$

We observe that since $\mathcal{B}_0$ is closed in $\mathcal{B}$, it is of second category and then $\mathcal{T}$ cannot be residual in $\mathcal{B}$.

\bigskip

In the context of Bloch functions,  it is natural to  study the Bergman functions. Let us introduce the notion of  Bergman functions. For $1\leq p <\infty$, $L^p(D,dA)$ will denote the Banach space of Lebesgue measurable functions on $D$, with the  $p$-norm 
$$\|f\|_p=\left(\int_D|f(z)|^pdA(z)\right)^p,$$ 
where $dA$ is the normalized area measure on $D$, that is, $dA(z)=\frac{1}{\pi}dxdy = \frac{1}{\pi}rdrd\theta,\,\,\,z=x+iy=re^{i\theta}.$ A function  $f:D\sobre\C$  is  a {\it Bergman function}  if $f$ is analytic  function in the open unit disk $D$ and $f \in L^p(D,dA).$ In other words, for $ 0 < p < \infty$, the Bergman space is the space of all analytic  functions $f$ in D for which the $p$-norm is finite. A  Bergman space is a function space of analytic  functions in the unit disk $D$ of the complex plane that are sufficiently well-behaved at the boundary that they are absolutely integrable.  We will denote   \emph{Bergman Space} by $A^p=\mathcal{H}(D)\cap L^p(D,dA),$ for  $ 1 \leq  p < \infty$.  So, the Bergman spaces are Banach spaces with $p$-norm.

It is known that the Bloch space is contained in all Bergman spaces, that is, $\mathcal{B}\subset A^p$, for all $0< p< \infty$ and the   inclusion is strict. Indeed, for a fixed $1\leq p <\infty$, consider $0<\beta<\dfrac{1}{p}$. Then the function $f_{(\beta)}:D\sobre\C$ defined by $f_{(\beta)}(z)=\dfrac{1}{(1-z)^{\beta}}$  is a Bergman function, but is not Bloch function. To see that $f_{(\beta)}\in A^p$, note that $\|f_{(\beta)}\|_p^p\leq\frac{1}{\pi}\int_0^{2\pi}\int_0^1\dfrac{r}{(1-r)^{\beta p}}drd\theta $, and this integral converges when $\beta<\frac{1}{p}$. To show that $f_{(\beta)}\notin \mathcal{B}$, note that for $x\in [0,1)$, $\lim_{x\tende 1^-}(1-x^2)|f_{\beta}'(x)|=+\infty$.

In the following, we will study the algebraic and topological  structures of the set   $\mathcal{W}^p=A^p\setminus \mathcal{B}$.

\begin{teo} For each $1\leq p<\infty$, the set $\mathcal{W}^p$ is spaceable. \end{teo}

\textbf{Proof:} To show that $\mathcal{W}^p=A^p\setminus \mathcal{B}$ is spaceable, we will show that the inclusion mapping $(\mathcal{B},\|\cdot\|_{\mathcal{B}})\hookrightarrow (A^p,\|\cdot\|_p)$ is continuous with unclosed range, and then apply \cite[Proposition 2.4]{kittim}. Let $f\in\mathcal{B},$  by \cite[Theorem 5.4]{zhu2}, for all $n\in\N$, the function $(1-|z|^2)^nf^{(n)}(z)$ is bounded on $D$, and then, by \cite[Theorem 4.28]{zhu2}, we have that $f\in A^p$. Let $g_f(z)=(1-|z|^2)f'(z)$. By the proof of \cite[Theorem 4.28]{zhu2}, there exists $C>0$ such that $\int_D|f|^pdA\leq C\int_D|g_f(w)|^p dA(w)$, and then we can find $M>0$ such that $\|f\|_p\leq M\|f\|_{\mathcal{B}}$, for all $f\in\mathcal{B}$. To show that $(\mathcal{B},\|\cdot\|_p)$ is not closed in $(A^p,\|\cdot\|_p)$, consider the function $f_{(\beta)}$ defined above ($0<\beta<\frac{1}{p}$). Let $(r_n)\subset [0,1]$ be a sequence such that $r_n\nearrow 1$, and consider the dilate of $f_{(\beta)}$ as $f_{n,(\beta)}(z)=\dfrac{1}{(1-r_n z)^{\beta}}$. It follows by \cite[Proposition 1.3]{hekozhu} that $f_{n,(\beta)}$ converges to $f_{(\beta)}$ in the $L^p$ norm, but $f_{(\beta)}\notin\mathcal{B}$. $\Box$

As a consequence we have that, for $1\leq p<\infty$, the set $\mathcal{W}^p$  is lineable, but the above proof does not display a linearly independent set. In the next proposition we  construct an infinite linearly independent set in $\mathcal{W}^p$. 

\begin{prop} For each $1\leq p < \infty$, there is an infinite linearly independent set $S$ in $A^p$ such that $[S]\subset\mathcal{W}^p\cup\{0\}$.\end{prop}

\textbf{Proof:} Let $S=\{f_{(\beta)}\,:\,0<\beta<\frac{1}{p}\}$, where $f_{(\beta)}$ is defined above. Let $f=\sum_{k=1}^{N}\lambda_k f_{(\beta_k)}$, where $\lambda_k\in\C$, $0<\beta_k<\frac{1}{p}$ and $\beta_k$ are distincts, for $k=1,\dots,n$. Suppose that $f(z)=0$, for all $z\in D$. Let $z_n=1-\dfrac{1}{n}$, for all $n\in\N$. Then $f(z_n)=\lambda_1n^{\beta_1}+\cdots+\lambda_Nn^{\beta_N}=0$, for all $n\in\N$.  Without loss of generality, we assume $\beta_1>\cdots>\beta_N$, and if $\lambda_1\neq 0$, then $n^{\beta_1}(\lambda_1+\lambda_2n^{\beta_2-\beta_1}+\cdots+\lambda_Nn^{\beta_N-\beta_1})\tende\infty$, as $n\tende\infty$, a contradiction. Then $\lambda_1=0$ and inductively we have that $\lambda_2=\cdots=\lambda_N=0$.

Consider now the vector space generated by $S$, denoted by $[S]$. It is clear that $[S]\subset A^p$. Let $f\in[S]$, that is, $f=\sum_{k=1}^{N}\lambda_k f_{(\beta_k)}$, where $\lambda_k\in\C$ and $0<\beta_k<\frac{1}{p}$, for $k=1,\dots,N$. Again,  we can assume that $\beta_1>\beta_2>\cdots>\beta_N$. We will show that $f\notin \mathcal{B}$. Note that 
$$(1-|z|^2)|f'(z)|\geq \dfrac{1-|z|^2}{|1-z|^{\beta_1+1}}\left(|\lambda_1\beta_1|-\dfrac{|\lambda_2\beta_2|}{|1-z|^{\beta_2-\beta_1}}-\cdots-\dfrac{|\lambda_N\beta_N|}{|1-z|^{\beta_N-\beta_1}}\right).$$ If in particular we take $z=x \in\R$, $|x|<1$, we see that $\sup_{z\in D}(1-|z|^2)|f'(z)|=+\infty$, and hence $f\notin \mathcal{B}$. $\Box$

We will finish  this section by proving that $\mathcal{W}^p$ is topologically big.

\begin{teo} For each $1\leq p<\infty$, the set $\mathcal{W}^p$ is residual in $A^p$. \end{teo}

\textbf{Proof:} Given an analytic function $f:D\sobre \C$, we will denote $g_f(z)=(1-|z|^2)f'(z)$. For each $n\in\N$, let $S_n=\{f\in A^p\,:\,\exists z\in D \textrm{ such that } |g_f(z)|>n\}$. Then it is clear that $A^p\setminus\mathcal{B}=\cap_{n=1}^{\infty}S_n$. If we show that each $S_n$ is open and dense in $A^p$, then it will imply that $A^p\setminus \mathcal{B}$ is residual in $A^p$.

If $f\in S_n$, then there exists $z_0\in D$ and $\delta>0$ such that $|g_f(z_0)|=n+\delta$. By \cite[Proposition 1.1]{hekozhu}, there exists $C>0$ such that $|h'(z_0)|\leq C\|h\|_p$, for all $h\in A^p$. So,  we choose $r$  such that $0<r<\dfrac{\delta}{C}$, then it follows that $B(f,r)\subset S_n$. Indeed, let $F\in B(f,r)$. Then $F-f=h$, with $\|h\|_p<r$. Now $|g_F(z_0)|\geq |g_f(z_0)|-|g_h(z_0)|>n+\delta-Cr>n$. 

To show that each $S_n$ is dense in $A^p$, let $f\in A^p$ and $\varepsilon>0$ be given. If $f\notin\mathcal{B}$, then it is clear that $B(f,\varepsilon)\cap S_n\neq\emptyset$. If $f\in\mathcal{B}$, let $0<\beta<\frac{1}{p}$ and $h=\dfrac{\varepsilon}{2}\dfrac{f_{(\beta)}}{\|f_{(\beta)}\|_p}$, $0<\beta<\frac{1}{p}$. Then there exists $z_0\in D$ such that $|g_h(z_0)|>n+\|f\|_{\beta}$. So if we take $F=f+h$, then $F\in B(f,\varepsilon)$ and $|g_F(z_0)|>n$. $\Box$

\medskip

We have studied the so-called unweighted Bergman spaces. The weighted Bergman spaces are defined as follows. For $\alpha>-1$, let $dA_{\alpha}(z)=(\alpha+1)(1-|z|^2)^{\alpha}dA(z)$. Then the \emph{weighted Bergman space} $A_{\alpha}^p$ is the set of all analytic functions $f:D\sobre\C$ such that $\int_D|f(z)|^pdA_{\alpha}(z)<\infty$. By the previous results we have the following corollary.

\begin{corol} For $\alpha>0$, the set $A_{\alpha}^p\setminus \mathcal{B}$ is spaceable and residual in $A_{\alpha}^p$.\end{corol}

\textbf{Proof:} Note that, for $\alpha>0$, $A^p\subset A_{\alpha}^p$, and the inclusion $(A^p,\|\cdot\|_p)\hookrightarrow (A_{\alpha}^p,\|\cdot\|_{p,\alpha})$. Now the arguments of the previous results apply analogously in this case. $\Box$

\section{Lorch Analytic Functions}

E. R. Lorch in \cite{lor} introduced  a definition of analytic mappings  (see Definition \ref{lanalytic}), that  have for their domains and ranges a complex commutative Banach algebra with identity. 

\begin{df}\label{lanalytic} Let $E$ be a complex  commutative Banach algebra  with identity. A mapping $f:E\sobre E$ \textbf{has a derivative in the sense of Lorch (an (L)-derivative) in} $\omega\in E$ if there exists $\zeta\in E$
such that $$\lim_{h\tende 0}\dfrac{\|f(\omega+h)-f(\omega)-\zeta\cdot h\|}{\|h\|}=0.$$
 We say that $f$ is \textbf{Lorch-analytic ((L)-analytic) in}  $E$,  if $f$ is \emph{(L)-analytic} in every point of $E.$
\end{df}

 We denote the space of all (L)-analytic mapping  from $E$ into $E$ by $\mathcal{H}_L(E, E)$.

Every (L)-analytic mapping is clearly  Fr\'echet differentiable, and hence analytic in the usual sense. The converse is not true,  because the fucntion  $F:\C^2\sobre\C^2$ given by $F(z,w)=(w,z)$, for all $(z,w)\in\C^2$, is such that $F \in  \mathcal{H}(E,E),$ but $F \notin \mathcal{H}_L(E,E)$ (see \cite{hilphi}).

The development of the primary aspects of the Lorch theory is  parallel that of the classical theory of analytic functions on complex variable. Lorch's definition allowed to consider several problems that were not studied in the scope of the standard theory of analytic mappings and several authors have  extended the theory.

In \cite{luizalex} and  \cite{luizalex2} the authors study  several topological and  algebraic properties of  the space $\mathcal{H}_L(E,E)$. Specially, they show that  $\mathcal{H}_L(E,E)$ is a closed subalgebra of the Fr\'echet algebra $\mathcal{H}_b(E,E),$  which consists of all entire functions of bounded type.  We study the existence of  large closed  linear space between these  two spaces.

 In \cite{lourvi2}, using the function $F$  defined above, we showed that $\mathcal{G}=\mathcal{H}(\C^2,\C^2)\setminus\mathcal{H}_L(\C^2,\C^2)$ is spaceable and not residual in $\mathcal{H}(\C^2,\C^2)$.
In this section, we want to investigate  the such spaces in case $E$ is   a general  infinite dimensional Banach algebra. We want to study two (different) sets: $\mathcal{G}(E)=\mathcal{H}(E,E)\setminus \mathcal{H}_L(E,E)$ and $\mathcal{G}_b(E)=\mathcal{H}_b(E,E)\setminus \mathcal{H}_L(E,E).$ We will show that both sets are nonempty, that $\mathcal{G}(E)$ is lineable and that $\mathcal{G}_b(E)$ is spaceable.

The following remarks will be useful for our results.

\begin{rema}\label{serie}
\begin{enumerate}

\item[(1)] \cite[Remark 2.3]{luizalex}. A holomorphic mapping $f:E\sobre E$ is (L)-analytic in $E$ if, and only if, there exists a unique sequence $(a_n)_{n\in\N}\subset E$ such that $\lim_{n\tende\infty}\|a_n\|^{\frac{1}{n}}=0$ and $f(z)=\sum_{n=0}^{\infty}a_nz^n$, for all $z\in E$.

\item[(2)] Let $\alpha\in\C$, $\alpha\neq 0$. If $f_{\alpha}(z)=f(\alpha z)$, for all $z\in E$, then $f\in \mathcal{H}_L(E,E)$ if and only if $f_{\alpha}\in \mathcal{H}_L(E,E)$.

\end{enumerate}
\end{rema}

\begin{teo}  Let $E$ be a complex  commutative Banach algebra with identity ${\rm e}$. The set $\mathcal{G}(E)$ is lineable. \end{teo}

\textbf{Proof.} Let $\omega\in \C$, $\omega\neq 0$. If we consider $a_n=\omega^n\,\textrm{e}.$  and $f(z)=\sum_{n=0}^{\infty}a_n z^n$, for all $z\in E$ we have that  $f\in\mathcal{G}(E)$. Consider $S=\{f_{\alpha}\,:\,\alpha\in\R\,\,\alpha\neq 0\}$. We will show that $S$ is linearly independent and $[S]\subset\mathcal{G}(E)\cup\{0\}$. 

Suppose that $\sum_{k=1}^N c_k f_{\alpha_k}=0$, with $c_k\in\C$, $\alpha_k\in\R$, $\alpha_k\neq 0$ for $k=1,\dots,N$. We can assume that $N\geq 2$ and $\alpha_1<\alpha_2<\cdots<\alpha_N$. Suppose that $c_N\neq 0$. Then it follows that 
$$\sum_{k=1}^N c_k f_{\alpha_k}(z)=\sum_{k=1}^N c_k \sum_{n=0}^{\infty}a_n\alpha_k^nz^n = 0,\textrm{ for all }z\in E.$$

Applying for $z=\lambda\textrm{e}$, where $|\lambda|<\delta$, for some $\delta>0$, we have that $$a_n(c_1\alpha_1^n+\cdots+c_N\alpha_N^n)=0,\textrm{ for all }n\in\N.$$

Since  $a_n\neq 0$, for all $n\in\N,$  we get that 
$$c_1\alpha_1^{n}+\cdots+c_N\alpha_N^{n}=0,\textrm{ for all }n\in\N.$$ 
As $c_N \neq 0,$ we have that
 $$\sum_{k=1}^{N-1}\dfrac{c_k}{c_N}\Big(\dfrac{\alpha_k}{\alpha_N}\Big)^{n}=-1,$$ By taking $n\tende\infty$ we find a contradiction, thus  $c_1=\cdots=c_N=0$ and  $S$ is linearly independent. 

Consider now $g=\sum_{k=1}^N c_k f_{\alpha_k}\in[S]$, with $c_k\in\C$, $\alpha_k\in\R$, $\alpha_k\neq 0$ for $k=1,\dots,N$, and suppose that $c_k\neq 0$, for $k=1,\dots,N$, and that $\alpha_1>\cdots>\alpha_N$. Then $$g(z)=\sum_{n=0}^{\infty}a_n(c_1\alpha_1^n+c_2\alpha_2^n+\cdots c_N\alpha_N^n)z^n.$$
If we call   $b_n=a_n(c_1\alpha_1^n+c_2\alpha_2^n+\cdots c_N\alpha_N^n), \forall n\in\N,$ to have  that $g\in\mathcal{G}(E)$, we must show that $\lim_{n\tende\infty}\|b_n\|^{\frac{1}{n}}\neq 0$.
 Observe that $\|b_n\|^{\frac{1}{n}}=\|a_n\|^{\frac{1}{n}}\|c_1\alpha_1^n+c_2\alpha_2^n+\cdots c_N\alpha_N^n\|^{\frac{1}{n}}.$ Then $$\|b_n\|^{\frac{1}{n}}=\|a_n\|^{\frac{1}{n}}|\alpha_1|\big\|c_1+c_2\big(\dfrac{\alpha_2}{\alpha_1}\big)^n+\cdots c_N\big(\dfrac{\alpha_N}{\alpha_1}\big)^n\big\|^{\frac{1}{n}}.$$ So, if $\lim_{n\tende\infty}\|b_n\|^{\frac{1}{n}}=0,$  this would imply that the sequence $$\Big(\big\|c_1+c_2\big(\dfrac{\alpha_2}{\alpha_1}\big)^n+\cdots c_N\big(\dfrac{\alpha_N}{\alpha_1}\big)^n\big\|\Big)_{n\in\N}$$ is bounded, which cannot happen unless $c_2=c_3=\cdots=c_N=0$.$\Box$

Note that the function $f$ defined above does not belong to $\mathcal{H}_b(E,E)$. In order to construct an element in $\mathcal{G}_b(E)$, we need to guarantee the existence of a linear functional $\varphi\in E'$ as in the following Lemma.

\begin{lemma}\label{functional} Let $E$ be a complex Banach algebra with identity $\textrm{e}$ and  dim($E$)$>1$. Then there exists $\varphi\in E'$, $\varphi(\textrm{e})=1$ and $\varphi(x_0)=0$ for some invertible element $x_0\in E$.\end{lemma} 

\textbf{Proof:} Let $U(E)$ denote the set of all invertible elements in $E$. Since $U(E)$ is an open set,  there exists $r>0$ such that $B(\textrm{e},r)\subset U(E)$. Let $u\in E$, $u\neq \alpha \textrm{e}$, for all $\alpha\in\C$ and $\|u\|=1$. Let  $x_0=\textrm{e}-\lambda u$, for some $0<|\lambda|<r$. Thus $x_0 \in U(E)$ and $x_0\neq\alpha \textrm{e}$, for all $\alpha\in\C$. Consider $M=[x_0]$, the closed vector subspace of $E$ generated by $x_0$. Since $\textrm{e}\notin M$, then there exists $\psi\in E'$ such that $\psi(\textrm{e})\neq 0$ and $\psi(x)=0$, for all $x\in M$. Now, the functional $\varphi=\dfrac{\psi}{\psi(\textrm{e})}$ has the desired properties. $\Box$

\begin{teo}Let $E$ be a complex Banach algebra with identity $\textrm{e}$. The set $\mathcal{G}_b(E)$ is spaceable.\end{teo}

\textbf{Proof:}  Firstly, we will  show that $\mathcal{G}_b(E)\neq\emptyset$. Let $(b_n)_{n\in\N}$ be a sequence in $E$ such that $\lim_{n\tende\infty}\|b_n\|^{\frac{1}{n}}=0$ and $b_n\neq 0$, for all $n\in\N.$  Using the   Lemma \ref{functional}, there exists  $\varphi\in E'$  with  $\varphi(\textrm{e})=1$ and $\varphi(x_0)=0$ for some invertible element $x_0\in E.$ Now, it is  possible to construct a sequence of $n$-homogenous polynomials im $E$.
 For each $n\in\N$,  we define  $P_n:E\sobre E$, $P_n(z)=b_n\varphi(z)^n$, for all $z\in E$. Thus $\lim_{n\tende\infty}\|P_n\|^{\frac{1}{n}}=0.$  If $g(z)=\sum_{n=0}^{\infty}P_n(z)$, for all $z\in E$, we have  that $g\in\mathcal{H}_b(E,E)$. If $g\in \mathcal{H}_L(E,E)$, then by Remark \ref{serie}, $P_n(z)=a_n z^n$, for some $(a_n)_{n\in\N}\subset E$, and for all $z\in E$. On the 
other  hand $P_n(x_0)=a_n x_0^{n}=b_n \varphi(x_0)^{n}=0,$ so $a_n=0$.Since  $P_n(\textrm{e})=a_n=b_n\varphi(\textrm{e})^{n}$ it follows  $b_n=0$, which is a contradiction. Therefore, $g \in \mathcal{G}_b(E)$

 Now, we will show that the quotient space $\mathcal{H}_b(E,E)/\mathcal{H}_L(E,E)$ is infinite dimensional, and  using the  \cite[Theorem 7.4.1 ]{kittim}, we obtain  that $\mathcal{G}_b(E)$ is spaceable.
Consider the set of classes $\mathcal{C}=\{\widehat{g_{\beta}}:\beta\in\R, 0<\beta< 1\}$, where $g$ was defined above, and $g_{\beta}(z)=g(\beta z)$, for all $z\in E$ We aim thar  the set $\mathcal{C}$ is linearly independent. Indeed, suppose that $\sum_{k=1}^N c_k\widehat{g_{\beta_k}}= \widehat{0}$, for $c_1,\dots,c_N\in\C$, and for  $k=1,\dots,N$,  all $\beta_k$  are distinct. If  we call  $h=\sum_{k=1}^N c_k g_{\beta_k}$, so  $h \in\mathcal{H}_L(E,E)$ and  if  we denote  by $d_n=c_1\beta_1^n+\cdots+c_N\beta_N^n$,  we have that
 $$h(z)=\sum_{n=0}^{\infty}d_n b_n\varphi(z)^n,\,\forall z\in E.$$

 Since $ h \in H_b(E,E)$, it follows $\lim_{n\tende\infty}\|d_n b_n\|^{\frac{1}{n}}=0$.  If there is $n_0\in\N$ such that  $d_{n_0}\neq 0$, then by the same arguments above we have that $h\notin\mathcal{H}_L(E,E)$, unless $h=0$. So, if $h(z)=0$, for all $z\in E$, then applying for $z=\lambda\textrm{e}$, where $|\lambda|<\delta$, for some $\delta>0$, we have that $d_nb_n=0$, for all $n\in\N$,it  implies $d_n = 0 \, \, \forall n$, it is  contradiction. Then $d_n = 0$  for all $n, $  that means. $c_1\beta_1^n+\cdots+c_N\beta_N^n=0$, for all $n\in\N$. Since we can take all the $\beta_k$ distincts, it follows that $c_1=\cdots=c_N=0$. The result follows. $\Box$.

\noindent\emph{Contact:}

\noindent M. Lilian Louren\c co, University of S\~ao Paulo, SP, Brazil, e-mail: mllouren@ime.usp.br

\noindent Daniela M. Vieira, University of S\~ao Paulo, SP, Brazil, e-mail: danim@ime.usp.br

\noindent R. do Mat\~ao, 1010 - Butant\~a, S\~ao Paulo - SP, Brazil, 05508-090

\end{document}